\documentclass{amsart}
\usepackage{amssymb}
 
\newtheorem{theorem}{Theorem} 
\newtheorem{lemma}[theorem]{Lemma}

\theoremstyle{definition}

\newtheorem{example}[theorem]{Example}

\theoremstyle{remark}

\DeclareSymbolFont{AMSb}{U}{msb}{m}{n}
\DeclareMathSymbol{\F}{\mathbin}{AMSb}{"46}
\DeclareMathSymbol{\N}{\mathbin}{AMSb}{"4E}
\DeclareMathSymbol{\Z}{\mathbin}{AMSb}{"5A}
\DeclareMathSymbol{\R}{\mathbin}{AMSb}{"52}
\DeclareMathSymbol{\C}{\mathbin}{AMSb}{"43}

\begin{document} \title[Hilbert's 17th Problem for matrices]
{An elementary and constructive solution \\ to 
Hilbert's 17th Problem for matrices}

\author{Christopher J. Hillar}
\address{Department of Mathematics, Texas A\&M University, College Station, TX 77843.}
\email{chillar@math.tamu.edu}
 \thanks{The first author is supported under an NSF
 Postdoctoral Research Fellowship.  This research was 
conducted during the \textit{Positive Polynomials and Optimization} workshop at
the Banff International Research Station, October 7--12 (2006), Banff, Canada.} 

\author{Jiawang Nie}
\address{Institute for Mathematics and its Applications, University of Minnesota,
Minneapolis, MN 55455.}
\email{njw@ima.umn.edu}

\subjclass{12D15, 03C64, 13L05, 14P05, 15A21, 15A54}

\keywords{Artin's theorem, Hilbert's 17th problem, sums of squares, positive semidefinite matrix,
real closed field}

\begin{abstract}
We give a short and elementary proof of a theorem of
Procesi, Schacher and (independently) Gondard, Ribenboim that generalizes a 
famous result of Artin.  Let $A$ be an $n \times n$ symmetric matrix with entries in the polynomial
ring $\mathbb R[x_1,\ldots,x_m]$.  The result is that if $A$ is postive semidefinite 
for all substitutions $(x_1,\ldots,x_m) \in \mathbb R^m$, then
$A$ can be expressed as a sum of squares of symmetric matrices with entries 
in $\mathbb R(x_1,\ldots,x_m)$.  Moreover, our proof is constructive and 
gives explicit representations modulo the scalar case.
\end{abstract} 

\maketitle 


We shall give an elementary proof of the following theorem.  Recall that
a real matrix is \textit{positive semidefinite} if it is symmetric with all nonnegative 
eigenvalues.

\begin{theorem}\label{mainthm}
Let $A$ be a symmetric matrix with entries in
$\mathbb R[x_1,\ldots,x_m]$.  If $A$ is postive semidefinite 
for all substitutions $(x_1,\ldots,x_m) \in \mathbb R^m$, then
$A$ can be expressed as a sum of squares of symmetric matrices with entries 
in $\mathbb R(x_1,\ldots,x_m)$.
\end{theorem}

This generalizes the following famous result of Artin on nonnegative polynomials;
it is the starting point for a large body of work relating positivity and algebra. 

\begin{theorem}[Artin]\label{artin}
If $f \in \mathbb R[x_1,\ldots,x_n]$ is nonnegative for all substitutions 
$(x_1,\ldots,x_n) \in \mathbb R^n$, then $f$ is a sum of squares of 
rational functions in $\mathbb R(x_1,\ldots,x_n)$.
\end{theorem}

Theorem \ref{mainthm} was originally proved in \cite{gondard} and 
(within a general framework) in \cite{procesi}, although a 
formulation involving elements 
in a number field was already considered in \cite{ciampi}.  
Like Artin's result, it guarantees algebraic certificates to (matrix) nonnegativity. 
However, the known proofs are nonconstructive, 
employing either model theory \cite{gondard} or ultraproducts \cite{procesi}.  
In contrast, we use only basic facts about real closed fields 
and linear algebra to give an explicit and elegant proof of Theorem \ref{mainthm}.  

Recall that a field $F$ is \textit{real} if 
$-1$ is not a sum of squares in $F$, and a \emph{real closed field} $R$ is a real
field such that any algebraic extension of $R$ that is
real must be equal to $R$.  Real closed fields have a unique ordering,
the nonnegative elements being the squares.
For instance,  $\mathbb R(x_1,\ldots,x_m)$ is 
a real field and $\mathbb R$ is real closed.  
A \textit{principal minor} of a matrix is
a determinant of a submatrix determined by the same row and column indices.
The set of symmetric matrices over $\mathbb R$ with all principal 
minors nonnegative coincides with the set of positive semidefinite matrices (see for example \cite[p. 405]{HJ1}), 
a fundamental relationship we exploit below.
We will prove the following generalization of Theorem \ref{mainthm} to the setting of real fields.



\begin{theorem}\label{genmainthm}
Let $F$ be a real field and let $A$ be a symmetric matrix with
entries in $F$.  If the principal minors of $A$ can be expressed as sums of 
squares in $F$, then $A$ is a sum of squares of symmetric matrices with 
entries in $F$.
\end{theorem}

To see see how Theorem \ref{mainthm} follows from Theorem \ref{genmainthm},
consider a principal minor $p(x_1,\ldots,x_m) \in \mathbb R[x_1,\ldots,x_m]$ of the
matrix $A$.  By assumption, it will be nonnegative for all substitutions 
$(x_1,\ldots,x_m) \in \mathbb R^m$, and therefore, Artin's theorem implies that
it is a sum of squares of rational functions.  We may now invoke Theorem \ref{genmainthm}.

As another application, consider positive semidefinite matrices $A \in \mathbb Q^{n \times n}$.  
Standard matrix theory allows one to
write $A = B^2$ for a symmetric
$B$ with entries that are algebraic numbers; however, Theorem \ref{genmainthm}
tells us that $A$ is actually a sum of squares of rational matrices.  This follows since any nonnegative 
rational number $a/b = ab/b^2$ can be written as a sum of four rational squares by
Lagrange's theorem.  

To prove Theorem \ref{genmainthm}, we begin with a lemma.  
For the basic theory of real closed fields (RCF) we will need,
we refer the reader to  \cite{Lang,Marker}.  
The main observation is that 
a symmetric matrix $A \in R^{n \times n}$ that has all nonnegative principal minors
is diagonalizable over $R$ with nonnegative eigenvalues, just as is the case for $\mathbb R$.

\begin{lemma}\label{minpolylem}
Suppose that $A$ satisfies the statement of Theorem \ref{genmainthm}.
Then the minimal polynomial $p(t) \in F[t]$ of $A$ 
is of the form: \[ p(t) = \sum_{i=0}^{m} (-1)^{m-i} a_i t^i = t^m - a_{m-1}t^{m-1} + \cdots + (-1)^m a_0\]
for $a_i$ that are sums of squares of elements of $F$.  Moreover, $a_1 \neq 0$.
\end{lemma}

\begin{proof}
Express the minimal polynomial of $A$ as in the statement of the theorem.
We first make the following observation.
Let $R$ be any real closure of $F$; this induces an ordering on 
$R$, in which the principal minors of $A$ are nonnegative (they
are sums of squares).
Since $A$ is diagonalizable over $R$ and has nonnegative eigenvalues, it follows that
each $a_i \geq 0$ and also that $p(t)$ has no repeated roots.

Suppose now that some $a_i$ was not a sum of squares in $F$. Then 
there is an ordering of $F$ with $a_i$ negative.  Let $R$ be a real closure of 
$F$ that extends the ordering on $F$.  By above, $a_i$ is nonnegative,
a contradiction.
To verify the second claim, first notice that $t^2$ does not divide $p(t)$ so that $a_0$ 
and $a_1$ cannot both be $0$.  In a real closure of $F$, 
the coefficient $a_1$ is a sum of products of (nonnegative) roots of 
$p(t)$.  It follows that if $a_1= 0$, we have $(-1)^ma_0 = p(0) = 0$.  Thus, $a_1 \neq 0$. 
\end{proof}

\begin{proof}[Proof of Theorem \ref{genmainthm}]
Let $A$ be a symmetric 
matrix satisfying the hypotheses of the theorem.  Also, let $p(t)$ be the 
minimal polynomial for $A$, which has the form prescribed by Lemma \ref{minpolylem}.
For notational simplicity, we assume that $m$ is odd, although the argument is
the same when $m$ is even. Since $p(A) = 0$, it follows that
\[(A^{m-1} + a_{m-2}A^{m-3} + \cdots + a_1I)A =  a_{m-1}A^{m-1} +
a_{m-3}A^{m-3}  + \cdots + a_0I.\]
Set $B = A^{m-1} +  \cdots + a_1I$, which is invertible (since $a_1 \neq 0$, in any
real closure of $F$, it is diagonalizable with strictly positive eigenvalues).  
Therefore, we have  
\begin{equation}\label{construction}
A = B \cdot \left(a_{m-1}B^{-2}A^{m-1} + a_{m-3}B^{-2} A^{m-3}  + \cdots + a_0B^{-2} \right).
\end{equation}
Since $B$ is a sum of squares and $B$ and $B^{-1}$ commute with $A$, the
result follows.
\end{proof}


Notice that our argument gives a commuting sum of squares representation,
the existence of which was also observed in \cite{procesi}.  
We close with two examples to illustrate the construction from our proof.

\begin{example}
The following symmetric matrix is always positive semidefinite:
\[
A = \begin{bmatrix}
1  & x_1x_2 \\ x_1x_2 & 1+x_1^4x_2^2+x_1^2x_2^4
\end{bmatrix}.
\]
However, it is not a sum of squares of matrix polynomials.
To see this, let $\textbf{x}  = [1,-1]^T$ and suppose that $A$ is a  sum
of polynomial squares; then so is the polynomial  
$f(x_1,x_2) = \textbf{x}^T A \textbf{x} = 2 +x_1^4x_2^2 +x_1^2x_2^4-2x_1x_2.$
Thus, we can express  $f  = \sum_{i=1}^n p_i^2$
for some polynomials $p_i$ with $\deg(p_i) \leq 3$. 
Comparing coefficients,
$p_i$ cannot contain the monomials $x_1^3$, $x_2^3$, $x_1^2$, $x_2^2$, $x_1x_2$, 
$x_1$, or $x_2$ so that we can write $p_i = a_i + b_i x_1^2x_2+c_i x_1x_2^2$
for some $a_i, b_i, c_i \in \mathbb R$. 
However, then we cannot produce the term $-2x_1x_2$ in $f$, a contradiction.
Similarly,  $\det(A)$ is not a sum of polynomial squares.  It
is, however, a sum of rational squares since $(x_1^2+x_2^2)\det (A)$ equals:
\begin{align*}
&  \left(x_2-\frac{1}{2}
x_1^2x_2\right)^2+\left(x_1-\frac{1}{2} x_1x_2^2\right)^2
+2\left(x_1x_2-\frac{1}{2}x_1x_2^3-\frac{1}{2}x_1^3x_2\right)^2 + \\
& \qquad \frac{3}{4}\left(x_1^2x_2^4+x_1^4x_2^2\right) +
\frac{1}{2}\left(x_1x_2^3+x_1^3x_2\right)^2.
\end{align*}
Since $A^2 - \text{tr}(A) A+ \det(A) I = 0$, we have the rational squares representation:
\[ A =  \text{tr}(A) \left[ \left(\text{tr}(A)^{-1}A \right)^2 +
  \det(A) \left(\text{tr}(A)^{-1} I \right)^2 \right]. \qed \]
\end{example}

\begin{example}
The following matrix is positive semidefinite for all substitutions:
\[ A = \begin{bmatrix}
x_1^2+2x_3^2 &-x_1x_2 & -x_1x_3 \\
-x_1x_2 & x_2^2+2x_1^2 & -x_2x_3 \\
-x_1x_3 & -x_2x_3 & x_3^2+2x_2^2
\end{bmatrix}, \]
but it is not a sum of polynomial squares \cite{choi}.
Its minimal polynomial has coefficients
\begin{align*}
a_2 = \ & 3x_3^2+3x_2^2+3x_1^2, \ 
a_1 = 2x_2^4+6x_1^2x_3^2+6x_1^2x_2^2+2x_1^4+2x_3^4+6x_2^2x_3^2, \\
a_0 = \ & 4x_1^4x_2^2+4x_3^2x_2^4+4x_3^4x_1^2+4x_3^2x_1^2x_2^2,
\end{align*}
which are all sums of squares. From formula (\ref{construction}), we have
\[ A = (A^2+a_1I)\left[
a_2 \left(A+a_1A^{-1}\right)^{-2} + a_0 \left(A^2+a_1I\right)^{-2}
\right].\qed\]
\end{example}

We would like to thank Konrad Schm\"udgen for discussing this problem with us.


\end{document}